\documentclass{my_fundam}

\usepackage{graphicx}
\usepackage{cite}
\usepackage{hyperref}

\begin{document}

\setcounter{page}{1}

\title{Non-differentiable Solutions for Local Fractional Nonlinear\\ 
Riccati Differential Equations}

\author{Xiao-Jun Yang\\
School of Mechanics and Civil Engineering,
China University of Mining and Technology\\
Xuzhou 221116, People's Republic of China\\
dyangxiaojun@163.com	
\and H. M. Srivastava\\
Department of Mathematics and Statistics, University of Victoria\\
Victoria, British Columbia V8W 3R4, Canada\\
{\it and} China Medical University, Taichung 40402, Taiwan, Republic of China\\
harimsri@math.uvic.ca
\and Delfim F. M. Torres\\
Center for Research and Development in Mathematics and Applications (CIDMA)\\
Department of Mathematics, University of Aveiro, 3810-193 Aveiro, Portugal\\
delfim@ua.pt
\and Yudong Zhang\\
School of Computer Science and Technology, Nanjing Normal University, 
Nanjing 210023\\
Department of Mathematics and Mechanics,
China University of Mining and Technology, Xuzhou 221008\\
People's Republic of China\\
zhangyudong@njnu.edu.cn}

\maketitle


\runninghead{X.-J. Yang et al.}{Non-differentiable Solutions
for Local Fractional Nonlinear Riccati Differential Equations}

\begin{abstract}
We investigate local fractional nonlinear Riccati differential equations (LFNRDE)
by transforming them into local fractional linear ordinary differential
equations. The case of LFNRDE with constant coefficients is considered
and non-differentiable solutions for special cases obtained.
\end{abstract}

\begin{keywords}
nonlinear Riccati equations,
non-differentiable functions,
local fractional derivatives.
\end{keywords}


\section{Introduction}

Ordinary differential equations (ODE) via local fractional derivatives \cite{1}
have played an important role in applied science, such as in fractal damped
vibrations \cite{2}, growth of populations \cite{3}, and fractal heat
transfers \cite{4}. The local fractional calculus is useful to describe
diffusion \cite{5,6,7,8}, heat conduction \cite{9},
waves \cite{10,11,12}, and other phenomena \cite{13}.
In \cite{14} the nonlinear Riccati differential equation (NRDE)
is written in the form
\begin{equation}
\label{eq1}
\frac{d\Lambda \left( \mu \right)}{d\mu }=\varpi _0 \left( \mu
\right)+\varpi _1 \left( \mu \right)\Lambda \left( \mu \right)+\varpi _2
\left( \mu \right)\Lambda ^2\left( \mu \right),
\end{equation}
where $\varpi_0 \left( \mu \right)\ne 0$ and $\varpi_2\left(\mu\right)\ne 0$.
The fractional nonlinear Riccati differential equation (FNRDE),
proposed in \cite{15}, takes the form
\begin{equation}
\label{eq2}
\frac{d^\upsilon \Theta \left( \mu \right)}{d\mu ^\upsilon }=\varpi _0
\left( \mu \right)+\varpi _1 \left( \mu \right)\Theta \left( \mu
\right)+\varpi _2 \left( \mu \right)\Theta ^2\left( \mu \right),
\end{equation}
where $\upsilon $ is the order of the fractional operator \cite{1,16,17,18},
$\varpi_0 \left( \mu \right)\ne 0$ and $\varpi _2 \left( \mu \right)\ne 0$.
Here we consider the local fractional nonlinear Riccati
differential equation (LFNRDE)
\begin{equation}
\label{eq3}
\frac{d^\zeta \Phi \left( \mu \right)}{d\mu ^\zeta }=\varpi _0 \left( \mu
\right)+\varpi _1 \left( \mu \right)\Phi \left( \mu \right)+\varpi _2 \left(
\mu \right)\Phi ^2\left( \mu \right),
\end{equation}
where $\zeta$ is the fractal dimension of the local fractional operator
\cite{1,2,3,4,5,6,7,8,9,10,11,12,13}, $\varpi _0 \left( \mu \right)\ne 0$
and $\varpi _2 \left(\mu\right)\ne 0$. Our main aim is to study 
non-differentiable solutions of LFNRDE. The highlights of this paper are: 
we prove that the local fractional nonlinear equation can be converted
into an equivalent linear equation with local fractional derivative,
and that its analytical solution can be found by using 
the local fractional Laplace transform.

The article is organized as follows. 
Section~\ref{sec:2} recalls the necessary tools
of local fractional differentiation.
In Section~\ref{sec:3}, transformation of the LFNRDE 
into a local fractional linear ODE is presented.
The LFNRDE with constant coefficients is studied in Section~\ref{sec:4}.
In Section~\ref{sec:5}, non-differentiable solutions for a LFNRDE are
discussed. Finally, the conclusion is outlined in Section~\ref{sec:6}.


\section{Mathematical tools}
\label{sec:2}

In this section, we introduce the basic theory of local fractional
differentiation. Suppose that a function $\Phi \left( \mu \right)$
is local fractional continuous in the domain $I=\left( {a,b} \right)$.
Then we write it as in \cite{1,2,3}:
\begin{equation}
\label{eq4}
\Phi \left( \mu \right)\in C_\zeta \left( {a,b} \right).
\end{equation}
Suppose that $\Phi \left( \mu \right)\in C_\zeta \left( {a,b} \right)$ and
$0<\zeta \le 1$. For $\delta >0$ and $0<\left| {\mu -\mu _0 } \right|<\delta$,
the limit
\begin{equation}
\label{eq5}
D^{\left( \zeta \right)}\Phi \left( {\mu _0 } \right)=\left.\frac{d^\zeta \Phi
\left( \mu \right)}{d\mu ^\zeta }\right|_{\mu =\mu _0} =\mathop
{\lim }\limits_{\mu \to \mu _0 } \frac{\Delta ^\zeta \left[ {\Phi \left( \mu
\right)-\Phi \left( {\mu _0 } \right)} \right]}{\left( {\mu -\mu _0 }
\right)^\zeta }
\end{equation}
exists and is finite \cite{1,2,3,5}, where $\Delta ^\zeta \left[{\Phi\left(\mu
\right)-\Phi \left({\mu _0}\right)}\right]\cong \Gamma \left( {1+\zeta }
\right)\left[ {\Phi \left( \mu \right)-\Phi \left( {\mu _0 } \right)}
\right]$. If $\mu \in \left( {a,b} \right)$, then
\begin{equation}
\label{eq6}
D^{\left( \zeta \right)}\Phi \left( \mu \right)=\frac{d^\zeta \Phi \left(
\mu \right)}{d\mu ^\zeta }=\Phi ^{\left( \zeta \right)}\left( \mu \right)
\end{equation}
(see \cite{1}). We say that $D_\zeta \left( {a,b} \right)$ is a $\zeta$
local fractional derivative set if (\ref{eq6}) is valid for any point
$\mu \in \left({a,b} \right)$.

\begin{theorem}[See \cite{1,2,3,5}]
Let $\Phi _1 \left( \mu \right)$, $\Phi _2 \left( \mu \right)
\in D_\zeta \left( {a,b} \right)$ be non-differentiable functions
defined on fractal sets. The following local fractional differentiation
rules hold:
\begin{description}
\item[(a)] $D^{\left( \zeta \right)}\left[ {\Phi _1 \left( \mu \right)\pm
\Phi_2\left( \mu \right)} \right]=D^{\left( \zeta \right)}\Phi _1\left(\mu
\right)\pm D^{\left( \zeta \right)}\Phi _2 \left( \mu \right)$;

\item[(b)] $D^{\left( \zeta \right)}\left[ {\Phi _1 \left( \mu \right)\Phi_2
\left( \mu \right)} \right]=\left[ {D^{\left( \zeta \right)}\Phi _1 \left(
\mu \right)} \right]\Phi _2 \left( \mu \right)+\Phi _1 \left( \mu
\right)\left[ {D^{\left( \zeta \right)}\Phi _2 \left( \mu \right)} \right]$;

\item[(c)] $D^{\left( \zeta \right)}\left[ \frac{\Phi _1 \left(
\mu \right)}{\Phi_2\left( \mu \right)} \right]=\frac{ {\left[
{D^{\left( \zeta \right)}\Phi _1\left( \mu \right)} \right]
\Phi _2 \left( \mu \right)-\Phi _1 \left( \mu\right)\left[
{D^{\left( \zeta \right)}\Phi _2 \left( \mu \right)} \right]}}{\Phi_2^2 \left(
\mu \right)}$, provided $\Phi _2 \left( \mu\right)\ne 0$.
\end{description}
\end{theorem}


\section{From LFNRDE to linear ODE via local fractional derivatives}
\label{sec:3}

Taking the function
\begin{equation}
\label{eq7}
\chi \left( \mu \right)=\Phi \left( \mu \right)\varpi _2 \left( \mu
\right),
\end{equation}
we find that
\begin{equation}
\label{eq8}
\begin{split}
D^{\left( \zeta \right)}\chi \left( \mu \right)&=D^{\left( \zeta
\right)}\left[ {\Phi \left( \mu \right)\varpi _2 \left( \mu \right)} \right]\\
&=\left[ {D^{\left( \zeta \right)}\Phi \left( \mu \right)}
\right]\varpi _2 \left( \mu \right)+\Phi \left( \mu \right)\left[ {D^{\left(
\zeta \right)}\varpi _2 \left( \mu \right)} \right]. \\
\end{split}
\end{equation}
Submitting Eq. (\ref{eq3}) to Eq. (\ref{eq8}), we have
\begin{equation}
\label{eq9}
\begin{split}
D^{\left( \zeta \right)}\chi \left( \mu \right)
&=\left[ {\varpi _0 \left(
\mu \right)+\varpi _1 \left( \mu \right)\Phi \left( \mu \right)+\varpi _2
\left( \mu \right)\Phi ^2\left( \mu \right)} \right]\varpi _2 \left( \mu
\right)+\Phi \left( \mu \right)\left[ {D^{\left( \zeta \right)}\varpi _2
\left( \mu \right)} \right] \\
&=\left[ {\varpi _0 \left( \mu \right)+\varpi _1 \left( \mu
\right)\Phi \left( \mu \right)+\varpi _2 \left( \mu \right)\Phi ^2\left( \mu
\right)} \right]\varpi _2 \left( \mu \right)+\left[ {\frac{D^{\left( \zeta
\right)}\varpi _2 \left( \mu \right)}{\varpi _2 \left( \mu \right)}}
\right]\chi \left( \mu \right).
\end{split}
\end{equation}
Taking
\begin{equation}
\label{eq10}
\Omega _1 \left( \mu \right)=\varpi _2 \left( \mu \right)\varpi _0 \left(
\mu \right)
\end{equation}
and
\begin{equation}
\label{eq11}
\Omega _2 \left( \mu \right)=\varpi _1 \left( \mu \right)+\frac{D^{\left(
\zeta \right)}\varpi _2 \left( \mu \right)}{\varpi _2 \left( \mu \right)},
\end{equation}
Eq. (\ref{eq9}) can be written as follows:
\begin{equation}
\label{eq12}
D^{\left( \zeta \right)}\chi \left( \mu \right)=\chi ^2\left( \mu
\right)+\Omega _2 \left( \mu \right)\chi \left( \mu \right)+\Omega _1 \left(
\mu \right).
\end{equation}
Set the function
\begin{equation}
\label{eq13}
\chi \left( \mu \right):=-\frac{D^{\left( \zeta \right)}\psi \left( \mu
\right)}{\psi \left( \mu \right)}.
\end{equation}
In the light of Eqs. (\ref{eq12})--(\ref{eq13}), we have
\begin{equation}
\label{eq14}
\begin{split}
D^{\left( \zeta \right)}\chi \left( \mu \right)
&=D^{\left( \zeta\right)}
\left[ {-\frac{D^{\left( \zeta \right)}\psi \left( \mu \right)}{\psi
\left( \mu \right)}} \right] \\
&=-\frac{D^{\left( {2\zeta } \right)}\psi \left( \mu \right)}{\psi
\left( \mu \right)}+\left[ {\frac{D^{\left( \zeta \right)}\psi \left( \mu
\right)}{\psi \left( \mu \right)}} \right]^2 \\
&=-\frac{D^{\left( {2\zeta } \right)}\psi \left( \mu \right)}{\psi
\left( \mu \right)}+\chi ^2\left( \mu \right)
\end{split}
\end{equation}
so that
\begin{equation}
\label{eq15}
\begin{split}
\frac{D^{\left( {2\zeta } \right)}\psi \left( \mu \right)}{\psi \left( \mu
\right)}&=-\Omega _2 \left( \mu \right)\chi \left( \mu \right)-\Omega _1
\left( \mu \right) \\
&=\Omega _2 \left( \mu \right)\frac{D^{\left( \zeta \right)}\psi
\left( \mu \right)}{\psi \left( \mu \right)}-\Omega _1 \left( \mu \right).
\end{split}
\end{equation}
Thus, we have
\begin{equation}
\label{eq16}
D^{\left( {2\zeta } \right)}\psi \left( \mu \right)-\Omega _2 \left( \mu
\right)D^{\left( \zeta \right)}\psi \left( \mu \right)+\Omega _1 \left( \mu
\right)\psi \left( \mu \right)=0.
\end{equation}
Note that Eq. (\ref{eq16}) is a local fractional linear ODE.
In view of Eqs. (\ref{eq7}) and (\ref{eq13}), we obtain
the following non-differentiable solution of Eq. (\ref{eq3}):
\begin{equation}
\label{eq17}
\Phi \left( \mu \right)=-\frac{D^{\left( \zeta \right)}\psi \left( \mu
\right)}{\varpi _2 \left( \mu \right)\psi \left( \mu \right)}.
\end{equation}


\section{LFNRDE with constant coefficients}
\label{sec:4}

We now consider the following LFNRDE with constant coefficients:
\begin{equation}
\label{eq18}
\frac{d^\zeta \Phi \left( \mu \right)}{d\mu ^\zeta }=\varpi _0 +\varpi _1
\Phi \left( \mu \right)+\varpi _2 \Phi ^2\left( \mu \right),
\end{equation}
where $\varpi _0 \ne 0$ and $\varpi _2 \ne 0$.
Following Eqs. (\ref{eq10}), (\ref{eq11}) and (\ref{eq16}), we have
\begin{equation}
\label{eq19}
\Omega _1 \left( \mu \right)=\varpi _2 \varpi _0
\end{equation}
and
\begin{equation}
\label{eq20}
\Omega _2 \left( \mu \right)=\varpi _1 ,
\end{equation}
such that
\begin{equation}
\label{eq21}
D^{\left( {2\zeta } \right)}\psi \left( \mu \right)-\varpi _1 D^{\left(
\zeta \right)}\psi \left( \mu \right)+\varpi _2 \varpi _0 \psi \left( \mu
\right)=0
\end{equation}
subject to the initial value conditions
\begin{equation}
\label{eq22}
D^{\left( \zeta \right)}\psi \left( 0 \right)=\alpha,
\quad \psi \left( 0
\right)=\beta,
\end{equation}
where
\begin{equation}
\label{eq23}
\Phi \left( \mu \right)=-\frac{D^{\left( \zeta \right)}\psi \left( \mu
\right)}{\varpi _2 \psi \left( \mu \right)}.
\end{equation}
Taking the local fractional Laplace transform (LFLT) \cite{1}
of (\ref{eq21}), we obtain
\[
\psi \left( s \right)=\frac{\alpha s^\zeta +\beta \left( {1+\varpi _1 }
\right)}{s^{2\zeta }-\varpi _1 s^\zeta +\varpi _2 \varpi _0 },
\]
where $s^\zeta $ is the local fractional Laplace operator.
When $\Sigma =\varpi _1^2 -4\varpi _2 \varpi _0 >0$, taking the inverse
LFLT, we obtain
\begin{equation}
\label{eq24}
\psi \left( \mu \right)=A_0 E_\zeta \left( {-C_0 \mu ^\zeta } \right)+B_0
E_\zeta \left( {-D_0 \mu ^\zeta } \right)
\end{equation}
such that
\begin{equation}
\label{eq25}
\Phi \left( \mu \right)=\frac{A_0 C_0 E_\zeta \left( {-C_0 \mu ^\zeta }
\right)+B_0 D_0 E_\zeta \left( {-D_0 \mu ^\zeta } \right)}{\varpi _2 \left[
{A_0 E_\zeta \left( {-C_0 \mu ^\zeta } \right)+B_0 E_\zeta \left( {-D_0 \mu
^\zeta } \right)} \right]},
\end{equation}
where
\begin{equation}
\label{eq26}
A_0 =\frac{\alpha }{2}+\frac{\beta \left( {1+\varpi _1 }
\right)-\frac{\alpha \varpi _1 }{2}}{\sqrt {\varpi _1^2 -4\varpi _2 \varpi
_0 } },
\end{equation}
\begin{equation}
\label{eq27}
B_0 =\frac{\alpha }{2}-\frac{\beta \left( {1+\varpi _1 }
\right)-\frac{\alpha \varpi _1 }{2}}{\sqrt {\varpi _1^2 -4\varpi _2 \varpi
_0 } },
\end{equation}
\begin{equation}
\label{eq28}
C_0 =\frac{\varpi _1 +\sqrt {\varpi _1^2 -4\varpi _2 \varpi _0 } }{2}
\end{equation}
and
\begin{equation}
\label{eq29}
D_0 =\frac{\varpi _1 -\sqrt {\varpi _1^2 -4\varpi _2 \varpi _0 } }{2}.
\end{equation}
For $\Sigma =\varpi _1^2 -4\varpi _2 \varpi _0 <0$, after taking the inverse
LFLT, we have
\begin{equation}
\label{eq30}
\psi \left( \mu \right)=A_1 E_\zeta \left( {-C_1 \mu ^\zeta } \right)+B_1
E_\zeta \left( {-D_1 \mu ^\zeta } \right)
\end{equation}
such that
\begin{equation}
\label{eq31}
\Phi \left( \mu \right)=-\frac{A_1 C_1 E_\zeta \left( {-C_1 \mu ^\zeta }
\right)+B_1 D_1 E_\zeta \left( {-D_1 \mu ^\zeta } \right)}{\varpi _2 \left[
{A_1 E_\zeta \left( {-C_1 \mu ^\zeta } \right)+B_1 E_\zeta \left( {-D_1 \mu
^\zeta } \right)} \right]},
\end{equation}
where
\begin{equation}
\label{eq32}
A_1 =\frac{\alpha }{2}+\frac{\beta \left( {1+\varpi _1 }
\right)-\frac{\alpha \varpi _1 }{2}}{i^\zeta \sqrt {4\varpi _2 \varpi _0
-\varpi _1^2 } },
\end{equation}
\begin{equation}
\label{eq33}
B_1 =\frac{\alpha }{2}-\frac{\beta \left( {1+\varpi _1 }
\right)-\frac{\alpha \varpi _1 }{2}}{i^\zeta \sqrt {4\varpi _2 \varpi _0
-\varpi _1^2 } },
\end{equation}
\begin{equation}
\label{eq34}
C_1 =\frac{\varpi _1 +i^\zeta \sqrt {4\varpi _2 \varpi _0 -\varpi _1^2 }
}{2}
\end{equation}
and
\begin{equation}
\label{eq35}
D_1 =\frac{\varpi _1 -i^\zeta \sqrt {4\varpi _2 \varpi _0 -\varpi _1^2}}{2}
\end{equation}
with the fractal imaginary $i^\zeta $.
By setting $\Sigma :=\varpi _1^2 -4\varpi _2 \varpi _0 =0$,
and taking the inverse LFLT, we obtain
\begin{equation}
\label{eq36}
\Phi \left( \mu \right)=\alpha E_\zeta \left( {-\frac{\varpi _1 }{2}\mu^\zeta }
\right)+\left( {\beta \left( {1+\varpi _1 } \right)-\frac{\varpi _1\alpha }{2}}
\right)\frac{\mu ^\zeta }{\Gamma \left( {1+\zeta }\right)}
E_\zeta \left( {-\frac{\varpi _1 }{2}\mu ^\zeta } \right).
\end{equation}


\section{An illustrative example}
\label{sec:5}

\begin{figure}
\centering
\includegraphics[width=3in]{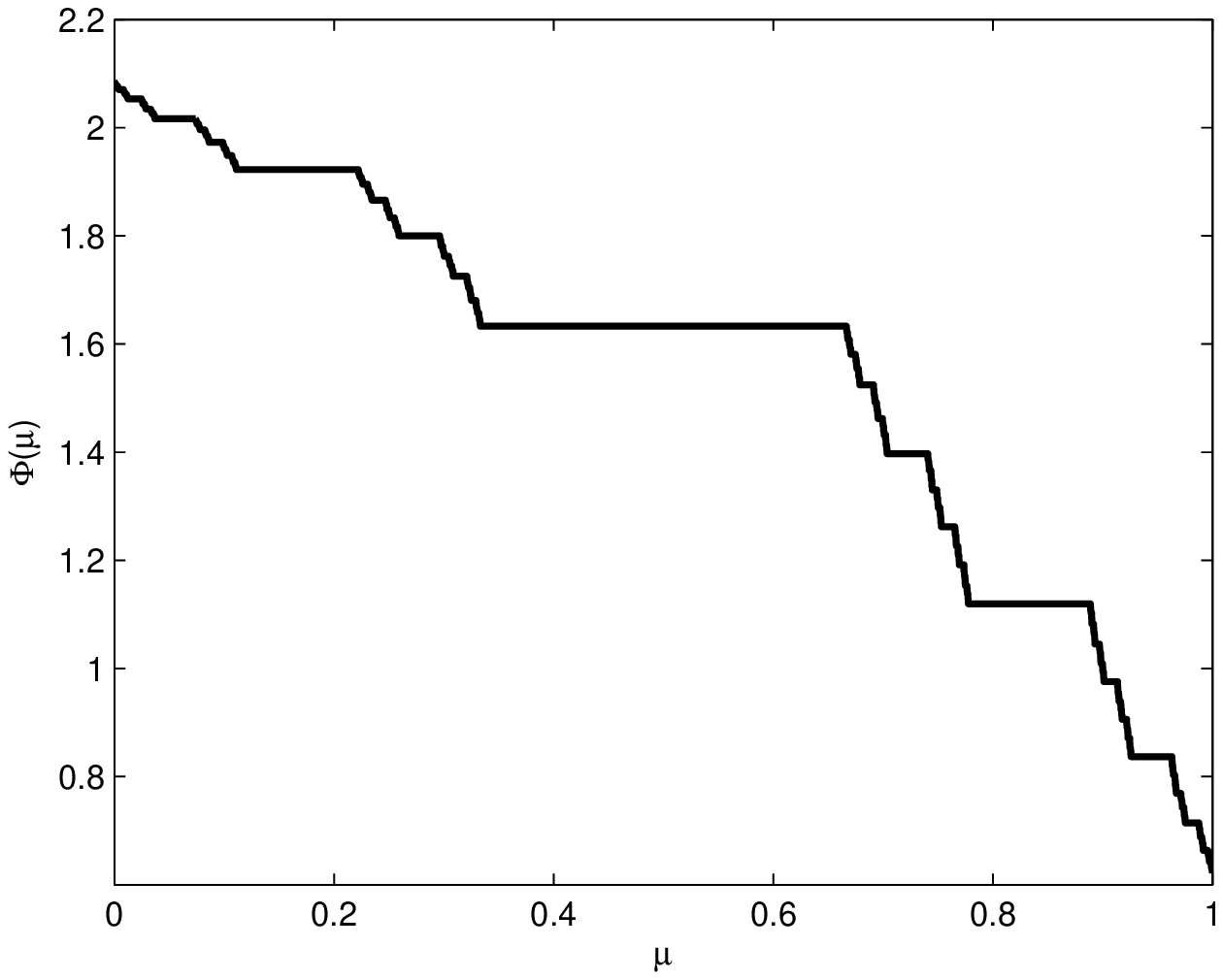}
\caption{Solution $\Phi \left( \mu \right)$ 
of (\ref{eq37})--(\ref{eq38}) 
when $\zeta =\ln 2/\ln 3$.}\label{figure1}
\end{figure}

Let us consider the LFNRDE with constant coefficients
\begin{equation}
\label{eq37}
\frac{d^\zeta \Phi \left( \mu \right)}{d\mu^\zeta}
=1+3\Phi \left(\mu\right)+\Phi ^2\left( \mu \right)
\end{equation}
subject to the initial condition
\begin{equation}
\label{eq38}
\Phi \left( 0 \right)=1.
\end{equation}
In light of Eqs. (\ref{eq21}) and (\ref{eq23}),
Eq. (\ref{eq38}) can be transformed into
\begin{equation}
\label{eq39}
D^{\left( {2\zeta } \right)}\psi \left( \mu \right)-3D^{\left( \zeta
\right)}\psi \left( \mu \right)+\psi \left( \mu \right)=0,
\end{equation}
where
\begin{equation}
\label{eq40}
D^{\left( \zeta \right)}\psi \left( 0 \right)=-\psi \left( 0 \right)=\xi.
\end{equation}
By setting $D^{\left( \zeta \right)}\psi \left( 0 \right)=-\psi \left( 0
\right)=\beta =-\alpha $, we have
\begin{equation}
\label{eq41}
\Sigma =\varpi _1^2 -4\varpi _2 \varpi _0=5,
\end{equation}
\begin{equation}
\label{eq42}
A_0 =\frac{\left( {\sqrt 5 -11} \right)\alpha }{2\sqrt 5 },
\end{equation}
\begin{equation}
\label{eq43}
B_0 =\frac{\left( {\sqrt 5 -5} \right)\alpha }{2\sqrt 5 },
\end{equation}
\begin{equation}
\label{eq44}
C_0 =\frac{3+\sqrt 5 }{2}
\end{equation}
and
\begin{equation}
\label{eq45}
D_0 =\frac{3-\sqrt 5 }{2}.
\end{equation}
Thus, making use of Eqs. (\ref{eq41})--(\ref{eq45}), we have
\[
\Phi \left( \mu \right)=\frac{\left( {\sqrt 5 -11} \right)\left( {3+\sqrt 5
} \right)E_\zeta \left( {-\frac{3+\sqrt 5 }{2}\mu ^\zeta } \right)+\left(
{\sqrt 5 -5} \right)\left( {3-\sqrt 5 } \right)E_\zeta \left(
{-\frac{3-\sqrt 5 }{2}\mu ^\zeta } \right)}{2\left( {\sqrt 5 -11}
\right)E_\zeta \left( {-\frac{3+\sqrt 5 }{2}\mu ^\zeta } \right)+2\left(
{\sqrt 5 -5} \right)E_\zeta \left( {-\frac{3-\sqrt 5 }{2}\mu ^\zeta }
\right)}
\]
and its corresponding graph, when $_{ }\zeta =\ln 2/\ln 3$,
is shown in Figure~\ref{figure1}.

\section{Conclusion}
\label{sec:6}

The nonlinear Riccati differential equation (NRDE) via local fractional
derivative is first proposed in this work. With the help of the properties
of the local fractional derivative, we can transform the LFNRDE into a local
fractional ODE. The non-differentiable solution of the LFNRDE with constant
coefficients is presented. An illustrative example with a Cantor-like
function is given. We claim that our method may be applied to improve
the performance in image classification, such as tea classification \cite{19,20},
pathological brain detection \cite{21}, fruit classification \cite{22,23}, 
and Alzheimer's disease detection \cite{24,25}. This is supported
by the fact that fractional nonlinear equations are commonly used 
in image processing. Our method can help to improve the solution 
accuracy, thus increasing the image classification.


\section*{Acknowledgements}

This research was partially supported
by the Open Project Program of the State Key Lab of CAD\&CG, 
Zhejiang University (A1616), and Natural Science Foundation 
of Jiangsu Province (BK20150983).
Torres was supported by CIDMA
and FCT within project UID/MAT/04106/2013.



\end{document}